\documentclass[]{interact}

\usepackage{epstopdf}
\usepackage[caption=false]{subfig}

\usepackage[numbers,sort&compress]{natbib}
\bibpunct[, ]{[}{]}{,}{n}{,}{,}
\renewcommand\bibfont{\fontsize{10}{12}\selectfont}

\theoremstyle{plain}
\newtheorem{theorem}{Theorem}[section]

\theoremstyle{definition}

\theoremstyle{remark}

\begin{document}

\title{Expected Time to Extinction of SIS Epidemic Model Using Quasy--
		Stationary Distribution}

\author{
\name{Kurnia Susvitasari\textsuperscript{1}\thanks{K.~Susvitasari. Email: susvitasari@sci.ui.ac.id}}
\affil{\textsuperscript{1} School of Mathematical Science, University of Nottingham, UK}}

\maketitle

\begin{abstract}

	We study that the breakdown of epidemic depends on some parameters, that is
	expressed in epidemic reproduction ratio number. It is noted that when $R_0 $
	exceeds 1, the stochastic model have two different results. But, eventually
	the extinction will be reached even though the major epidemic occurs. The question
	is how long this process will reach extinction. In this paper, we will focus on
	the Markovian process of SIS model when major epidemic occurs. Using
	the approximation of quasi--stationary distribution, the expected mean time
	of extinction only occurs when the process is one step away from being extinct.
	Combining the theorm from Ethier and Kurtz, we use CLT to find the approximation
	of this quasi distribution and successfully determine the asymptotic mean
	time to extinction of SIS model without demography.

\end{abstract}


\section{Introduction}

	According to the threshold theorem, epidemic can only occur if the initial number
	of susceptibles is larger than some critical values, which depends on the 
	parameters under cosideration \cite{ball} .
	Usually, it is expressed in terms of epidemic reproductive ratio number, $R_0$.
	
	Susvitasari \cite{kurnia}
	showed that both deterministic and stochastic models performed similar results
	when $R_0 \leq 1$, i.e. the disease--free stage in the epidemic. But then, when
	$R_0 > 1$, the deterministic and stochastic models had different interpretations.
	In the deterministic model, both SIS and SIR showed an outbreak of the 
	disease, and after some time $t$, the disease persisted and reached equilibrium
	stage, i.e. endemic. The stochastic model, on the other hand, had different
	interpretation. There are two possible outcomes of this approach.
	First, the infection may die out in the first cycle. If it did, it would happen very
	quickly since the time of the disease removed must less than the time of
	infectee--susceptible contact. Second, if it survives the first cycle, the outbreak
	was likely to occur, but after some time $t$, it would reach equilibrium as the
	deterministic version. In fact, the stochastic model would mimic the deterministic
	path and scattered randomly around its equilibrium. Furthermore, by letting
	population size be large and ignoring the initial value of infectious individuals,
	as $t \rightarrow \infty$, the empirical distribution of $n^{-1} Y_n(t)$, where
	$Y_n(t)$ is the number of infectious individuals at time $t$ with $n$ population
	size followed normal distribution in endemic phase as seen in Figure 
	\ref{fig:sis.empirical}.
	
	In this paper, we will focus on the endemic stage, the stage where $R_0 > 1$
	and the outbreak occurs in SIS epidemic model without demography to
	determine the expected mean time to extinction of this model.

\section{The Birth and Death Process in SIS}

\begin{figure}[ht]
		\centering
		\resizebox*{10cm}{!}{\includegraphics{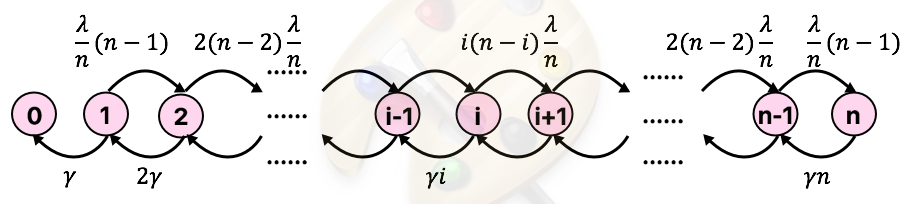}}
		\caption{Birth--death of SIS without demography transition process}
		\label{fig: sis birth--death}
	\end{figure}

	Suppose that $\lbrace Y_n(t): t \geq 0 \rbrace$ denotes the Markovian process
	that represents the number of infected individuals in the population at time $t$.
	Then we define $S_n = \lbrace 0, 1, \ldots, n \rbrace$ as the state space of the 
	process, where state zero is an absorbing state. Figure \ref{fig: sis birth--death}
	is the illustration of this process.
	
	Consider the case where $R_0 > 1$ and the outbreak occurs. Eventually,
	given a necessary large $n$ and some time $t$, the process will not visit 
	state zero. So, let us redefine the birth--death process before by only
	considering the transient states. Let $\lbrace \hat{Y}_n(t): t \geq 0 \rbrace$
	be the modified process of $\lbrace Y_n(t): t \geq 0 \rbrace$ with state space
	$S_n' = S_n \setminus \lbrace 0 \rbrace$. The intuitive explaination is that 
	whenever the SIS process goes extinct, it immediately will restart again in state 1.
	It is not hard to show that the process $\lbrace \hat{Y}_n(t): t \geq 0 \rbrace$
	is time--reversible Markov chain. Therefore, there exist colections of equilibrium
	distributions $\pi_j$ for all $j \in S_n'$ that sums to unity and satisfy the 
	detailed--balanced conditions:
		\begin{enumerate}
			\item $(n-1) \displaystyle \frac{\lambda}{n} \hat{\pi}_i = 2 \gamma \hat{\pi}_2$,
			\item $\left( i(n-i) \displaystyle \frac{\lambda}{n} + \gamma i \right) 
					\hat{\pi}_i = (i-1)(n-i+1) \displaystyle \frac{\lambda}{n} \hat{\pi}_{i-1} + (i+1) 
					\gamma \hat{\pi}_{i+1}$ \qquad $i = 2,3, \ldots, n-1$,
			\item $n \gamma \hat{\pi}_n = (n-1) \displaystyle \frac{\lambda}{n} 
					\hat{\pi}_{n-1}$,
		\end{enumerate}
	with constraint $\displaystyle \sum_{i \in S_n'} \hat{\pi}_i = 1$. 
	Thus the solution is
		\begin{equation}
			\hat{\pi}_i = \frac{1}{i} \left( \frac{\lambda}{n \gamma} \right) ^{i-1} \frac{(n-1)!}{(n-i)!} 	
			\hat{\pi}_1 \qquad i = 1,2, \ldots, n.
		\end{equation}
	Therefore,
		\begin{equation}
			\hat{\pi}_1 = \left( \sum_{i = 1}^n \frac{(n-1)!}{i(n-i)!} \left( \frac{\lambda}
			{n \gamma} \right) ^{i-1} \right) ^{-1}.
		\end{equation}
		
	Now, let $T_{ext}^{(n)}$ denote the time to extinction in SIS epidemic model. 
	Recall that the process $\lbrace \hat{Y}_n(t): t \geq 0 \rbrace$ is time-reversible 
	on state $S_n'$. Then, the time to extinction $T_{ext}^{(n)}$ is almost surely 
	finite since the state $S_n'$ is visited finitely often.
	
	Suppose that $\lbrace \hat{Y}_n(t): t \geq 0 \rbrace$ is in equilibrium, then
	the mean time to extinction has intensity $\hat{\pi}_1$ times the rate of the 
	process getting absorbed to state $\lbrace 0 \rbrace$.
		\begin{align}
			E \left( T_{ext}^{(n)} \right) &= E( \text{interarival time of jump } 1 \rightarrow 0) 
														\nonumber \\
													&= \frac{1}{\gamma \hat{\pi}_1} = 
													\frac{1}{\gamma} \sum_{i = 1}^n \frac{(n-1)!}{i(n-i)!} 
													\left( \frac{\lambda}{n \gamma} \right) ^{i-1}.
		\end{align}
	As $n \rightarrow \infty, E \left( T_{ext}^{(n)} \right) \rightarrow
	\displaystyle -\frac{\log {(1- \lambda)}}{\lambda}$ when we set
	$R_0 < 1$. Ball et al. \cite{ball-britton-neal} showed the approximation of mean time extinction
	when $R_0 < 1, R_0 = 1,$ and $R_0>1$.
	
	In this paper, we will focus on the time to extinction when the
	process reaches endemic--stage, i.e. when $R_0 > 1$ and
	the process survives the first cycle of extinction using quasi--
	statioary distribution and property of central limit theorm.
	
\section{Expected Time of Extinction of SIS in Endemic--Stage}

	Before we start this section, we will introduce first the quasi--stationary 
	distribution.

\subsection{Quasi--Stationary Distribution}

	Recall the process $\lbrace Y_n(t): t \geq 0 \rbrace$ with finite state space
	$S_n$, where 0 is an absorbing state. For $j \in T = S_n \setminus 
	\lbrace 0 \rbrace$, let
		\begin{align}
			q_{i,j}(t) &= P(Y_n(t)=j \mid Y_n(t)>0, Y_n(0) = i) 
			\nonumber \\
			&= \frac{P_{i,j}(t)}{1-P_{i,0}(t)} \qquad i = 1,2, \ldots, n. \label{eq:prob.q_ij}
		\end{align}
	Since we concerns on the distribution of the process in endemic stage, we
	wish to find $\tilde{q}_j = \displaystyle \lim_{t \rightarrow \infty} q_{i,j}(t)$.
	
	We define $A$ and $T$ as sets of absorbing and transient state. Suppose
		$$\boldsymbol{Q} = \left[
			\begin{array}{cc}
				\boldsymbol{Q}_{AA} & \boldsymbol{0} \\
				\boldsymbol{Q}_{TA} & \boldsymbol{Q}_{TT}
			\end{array} \right] $$ 
	where $\boldsymbol{Q}$ is transition rate matrix, defined as $\boldsymbol{Q}
	= (q_{ij})_{i,j \in S_n}$. Let $\boldsymbol{W} = diag_{n \times n}(\boldsymbol{
	\hat{\pi}})$. Since $\boldsymbol{W}^{1/2} \boldsymbol{Q}_{TT} \boldsymbol{W}^
	{-1/2}$ is symmetry, according to spectral decomposition theorem,
		\begin{equation}
			\boldsymbol{W}^{1/2} \boldsymbol{Q}_{TT} \boldsymbol{W}^{-1/2} = 
				\sum_{i=1}^n \lambda_i \boldsymbol{u}_i \boldsymbol{u}_i^T, \nonumber
		\end{equation}
	where $\lambda_1, \lambda_2, \ldots, \lambda_n$ are eigenvalues of 
	$\boldsymbol{W}^{1/2} \boldsymbol{Q}_{TT} \boldsymbol{W}^{-1/2}$ and 
	$\boldsymbol{u}_1, \boldsymbol{u}_2, \ldots, \boldsymbol{u}_n$ are corresponding 
	orthonormal vectors of right-eigenvectors. Therefore,
		\begin{equation}
			\boldsymbol{Q}_{TT} = \sum_{i=1}^n \lambda_i \boldsymbol{E}_i
			\Rightarrow \boldsymbol{Q}_{TT} \cdot t = \sum_{i=1}^n \lambda_i t
			\cdot \boldsymbol{E}_i \Rightarrow exp(\boldsymbol{Q}_{TT}  t) =
			exp \left( \sum_{i=1}^n \lambda_i t \cdot \boldsymbol{E}_i \right),
		\end{equation}
	where $\bold{E}_i = \left( \bold{W}^{-1/2} \bold{u}_i \right) \left( \bold{u}_i^T 
	\bold{W}^{1/2} \right)$ and we can see in Appendix .........
	that $exp(\boldsymbol{Q}_{TT}  t)  =\displaystyle  \sum_{i=1}^n exp(\lambda_it) 
	\boldsymbol{E}_i$. Recall that forward Kolmogorov equation, $\bold{P}(t)$ must 
	satisfy $\bold{P}'(t) = \bold{P}(t) \bold{Q}$, where $\bold{P}(0) = \bold{I}_{n \times n}$
	with solution $\bold{P}(t) = exp(\bold{Q}t)$. So, for $i, j \in T$,
		\begin{align}
			P_{i,j}(t) &= P \lbrace Y_n(t) = j \mid Y_n(0) = i \rbrace \nonumber \\
			&= \left[ e^{\boldsymbol{Q}_{TT}t} \right] _{ij} 
			= \left[ \sum_{k=1}^n e^{\lambda_k t} \boldsymbol{E}_k \right] _{ij}.
			\label{eq:prob(t)}
		\end{align}
	
	According to Darroch and Seneta \cite{darroch-seneta}
	there exists a simple, real eigenvalue $\lambda_1 < 0$ such that all other 
	eigenvalues of $\boldsymbol{W}^{1/2} \boldsymbol{Q}_{TT} \boldsymbol{W}^{-1/2}$ 
	have real part less than $\lambda_1$. Therefore, eq. (\ref{eq:prob.q_ij}) becomes
		\begin{align}
			q_{i,j}(t) &= \frac{\displaystyle e^{\lambda_1 t} (\boldsymbol{E}_1)_{ij} + 
			e^{\lambda_2 t} (\boldsymbol{E}_2)_{ij} + \ldots + e^{\lambda_n t} 
			(\boldsymbol{E}_n)_{ij}}{\displaystyle \sum_{l=1}^n \lbrace e^{\lambda_1 t} 
			(\boldsymbol{E}_1)_{il} + e^{\lambda_2 t} (\boldsymbol{E}_2)_{il} + \ldots + 
			e^{\lambda_n t} (\boldsymbol{E}_n)_{il} \rbrace} \nonumber \\
			&= \frac{\displaystyle (\boldsymbol{E}_1)_{ij} + e^{(\lambda_2 - \lambda_1) t} 
			(\boldsymbol{E}_2)_{ij} + \ldots + e^{(\lambda_n - \lambda_1) t} 
			(\boldsymbol{E}_n)_{ij}}{\displaystyle \sum_{l=1}^n \lbrace (\boldsymbol{E}_1)_{il} + 
			e^{(\lambda_2 - \lambda_1) t} (\boldsymbol{E}_2)_{il} + \ldots + e^{(\lambda_n - 
			\lambda_1) t} (\boldsymbol{E}_n)_{il} \rbrace}.
		\end{align}	
	Consequently, as $t \rightarrow \infty, \tilde{q}_j = \displaystyle  \frac{(\bold{E}_1)_{ij}}
	{\sum_{l=1}^n (\bold{E}_1)_{il}}$.
	
	Note that $\boldsymbol{E}_1 = \left( \boldsymbol{W}^{-1/2} \boldsymbol{u}_1 \right) \left( 
	\boldsymbol{u}_1^T \boldsymbol{W}^{1/2} \right)$ and it turns out that $\left( 
	\boldsymbol{u}_1^T \boldsymbol{W}^{1/2} \right)$ is left eigenvector of 
	$\boldsymbol{Q}_{TT}$ (see Appendix B). In the stationary 
	scenario, when $t \rightarrow \infty$, each row of probability transition matrix 
	$\boldsymbol{P}(t)$ has nearly similar entries. So, it follows that the quasi-distribution 
	$\tilde{\boldsymbol{q}}$ is given by the left-eigenvector of $\boldsymbol{Q}_{TT}$ 
	corresponding to eigenvalue $\lambda_1$ with constraint $\tilde{\boldsymbol{q}} 
	\boldsymbol{1} = \displaystyle \sum_{k=1}^n \tilde{q}_k = 1$.

\subsection{Expected Time to Extinction by Quasi--Stationary Distribution}

	Suppose that $\left( Y_n(t) \mid Y_n(t) > 0, Y_n(0) = i \sim \tilde{ \boldsymbol{q}} 
	\right)$ and let $T_Q = \inf{ \lbrace t \geq 0 ; Y_n(t) = 0, Y_n(t) \sim 
	\tilde{\boldsymbol{q}}} \rbrace$. Since two events (infection and recovery) in the 
	epidemic model follows the Poisson process, there is zero probability that two or 
	more events occur at the same time. So, the extinction only occurs when the 
	process is one step away from being extinct. Therefore,	
		\begin{align}
		P(\text{not extinct at time }t) &= P \left( T_Q > t \right) \nonumber \\
		&= \sum_{j=1}^n \lim_{t \rightarrow \infty} P(Y_n(t) = j \mid Y_n(0) = i) \nonumber \\
		&= \sum_{j=1}^n \tilde{q}_j \cdot P(Y_n(t) > 0) 	
		\nonumber \\
		&= \sum_{j=1}^n \tilde{q}_j \cdot e^{\boldsymbol{Q}_{TT}t} 
		= \tilde{\boldsymbol{q}} \cdot e^{\boldsymbol{Q}_{TT}t} \cdot \boldsymbol{1}. 	
		\label{eq:time.extinct}
		\end{align}
	Recall that $\tilde{\boldsymbol{q}}$ is the left-eigenvector of $\boldsymbol{Q}_{TT}$, 
	corresponding to the eigenvalue $\lambda_1$. So, $$\tilde{\boldsymbol{q}} \cdot 
	\boldsymbol{Q}_{TT} = \lambda_1 \tilde{\boldsymbol{q}} \quad \Rightarrow \quad 
	\tilde{\boldsymbol{q}} \cdot \boldsymbol{Q}_{TT}^k = \lambda_1^k \tilde{\boldsymbol{q}} 
	\quad \Rightarrow \quad \tilde{\boldsymbol{q}} \cdot e^{\boldsymbol{Q}_{TT}t} = 
	e^{\lambda_1 t} \tilde{\boldsymbol{q}}.$$
	Hence, the equation (\ref{eq:time.extinct}) becomes
		\begin{align}
			P(\text{not extinct at time }t) &
			= e^{\lambda_1 t}. \label{eq:time.extinct.tidy}
		\end{align}
		
	Since the process is Poisson process, then the inter arrival time must follow
	exponential distribution. Therefore, in eq. (\ref{eq:time.extinct.tidy}), $T_Q$
	follows exponential distribustion with intensity $-\lambda_1$. But, recall the 
	memoryless property of the Poisson process and that the extinction only occurs 
	if the process is one step away from being extinct. Intuitively, it also means that 
	the extinction happens with intensity $\tilde{q}_1$ times the rate of the process 
	to jump into absorption state $0$, which is $\gamma$. On the other hand, 
	$\lambda_1 = -\gamma \tilde{q}_1$ (note that $\lambda_1 < 0$). 
	
	\begin{proof}
		Using the fact that $\tilde{\boldsymbol{q}}$ is the left-eigenvector of 
		$\boldsymbol{Q}_{TT}$,
		\[
			\tilde{\boldsymbol{q}} \cdot \boldsymbol{Q}_{TT} = \lambda_1 
			\tilde{\boldsymbol{q}} \Rightarrow
			\tilde{\boldsymbol{q}} \cdot \underbrace{ \boldsymbol{Q}_{TT} \cdot 	
			\boldsymbol{1}}_{(-\gamma,0,\ldots,0)^T} = \lambda_1 \underbrace{ 
			(\tilde{\boldsymbol{q}}  \boldsymbol{1})}_{=1} \Rightarrow
 			-\gamma \cdot \tilde{q}_1 = \lambda_1.
		\]
	\end{proof}
	Therefore, the mean time to extinction of the SIS process is
		\begin{align}
			E \left( T_Q \right) 
			&= \frac{1}{\gamma \tilde{q}_1}. \label{eq:sis.time.quasi}
		\end{align}
	where $\tilde{q}_1 = (\tilde{\boldsymbol{q}})_{1}$.

	The quasi--stationary is a very powerful approximation when we let population
	size large enough. Unfortunately, consider that by using this method, we need
	to capture all the possible transition rate for matrix $\bold{Q}$. Considering
	the large $n$, the size of matrix $\bold{Q}$ will also big and leads to inefficient 
	computation. Therefore, we should approximate the distribution of $\tilde{q}_1$.
	
	Now, suppose that we only concern on the size of infectious individuals in
	population when the outbreak occurs, i.e. $\lbrace Y_n(t) | Y_n(t) > 0, Y_n(t)
	\sim \bold{\tilde{q}} \rbrace$. By ignoring the initial value of infectives, the empirical
	distribution of $n^{-1} Y_n(t)$ follows normal distribution with mean converges
	to equilibrium point in deterministic model as seen in Figure 
	\ref{fig:sis.empirical}. It turns
	out that by applying Central Limit Theorem, we could approximate the quantity
	of $\tilde{q}_1$.

\section{Central Limit Theorem}
	
	We define $\lbrace \boldsymbol{X}_n(t); t \geq 0 \rbrace$ to be the CTMC on 
	$\mathbb{Z}^d$ with finite state space and finite number of possible transition 
	$\boldsymbol{l} \in \Delta \subseteq \mathbb{Z}^d$ in such a way that 
	$\beta_{\boldsymbol{l}}: 
	\mathbb{Z}^d \rightarrow [0, \infty), \sum_{l \in \Delta} 
	\beta_l(\boldsymbol{x}) < \infty$ for all $\boldsymbol{x} \in \mathbb{Z}^d$.
	Then let $F(\boldsymbol{x}) = \sum_{\boldsymbol{l} \in \Delta} 
	\boldsymbol{l} \beta_{\boldsymbol{l}} (\boldsymbol{x})$ as the drift function of 
	the process $\lbrace \boldsymbol{X}_n(t): t \geq 0 \rbrace$.
	
	Suppose that $F(\boldsymbol{x})$ is a continuous function and differentiable in 
	$\mathbb{Z}^d$ lattice such that $\partial F(\boldsymbol{x})$ exists and is 
	continuous and let $G(\boldsymbol{x})$ be $d \times d$ matrix function defined as
	$$G(\boldsymbol{x}) = \displaystyle \sum_{\boldsymbol{l} \in \Delta} \boldsymbol{l}^T 
	\boldsymbol{l} \beta_{\boldsymbol{l}} (\boldsymbol{x}).$$
	Further, let $\Phi(t,s)$, where $0 \leq s \leq t$, be the solution of the matrix 
	differential equation
		\begin{equation}
			\frac{\partial}{\partial t} \Phi(t,s) = \partial F(\boldsymbol{x}(t)) \cdot \Phi(t,s) 	
			\label{eq:clt.diff.phi}
		\end{equation}
	where $\Phi(s,s) = \boldsymbol{I}$.
	
	\begin{theorem} \label{thm:clt}
		(Ethier and Kurtz) Suppose that
		\begin{enumerate}
			\item $\displaystyle \lim_{n \rightarrow \infty} \sqrt{n} \left( \frac{\displaystyle 	
					\boldsymbol{X}_n(0)}{n} - \boldsymbol{x}_0 \right) = \boldsymbol{v}_0$ 
					is a constant, where $\boldsymbol{x}(t) = \boldsymbol{x}_0 + \displaystyle 
					\int_0^t F(\boldsymbol{x}(u)) du$,
			\item $\beta_{\boldsymbol{l}}$ and $\partial F(\boldsymbol{x})$ are continuous 
					functions.
		\end{enumerate}
		Then, for $t \geq 0$,
			\begin{equation}
				 \sqrt{n} \left( \frac{\boldsymbol{X}_n(t)}{n} - \boldsymbol{x}(t) \right) 	
				 \xrightarrow[\text{weakly}]{\quad D \quad} \boldsymbol{V}(t) \nonumber
			\end{equation}
		as $n \rightarrow \infty$, where $\lbrace \boldsymbol{V}(t): t \geq 0 \rbrace$ is a 
		zero--mean Gaussian process with covariance function given as
			\begin{equation}
				\text{cov}(\boldsymbol{V}(t), \boldsymbol{V}(s)) = \int_0^{\min{(t,s)}} 
				\Phi(t,u) \cdot G(\boldsymbol{x}(u)) \cdot \Phi (s,u)^T du. 
				\label{eq:clt.cov}
			\end{equation}
		The functions $\Phi(t,s)$ and $G(\boldsymbol{x})$ are as defined previously.
	\end{theorem}
	
	The proof of this theorem can be seen in Britton and Andersson
	\cite{britton-andersson}  or in
	Ethier and Kurtz \cite{ethier-kurtz}.
	
	Note that the covarian function in equation (\ref{eq:clt.cov}) is time-dependent. 
	It implies that the variance function is also time-dependent. Using equation 
	(\ref{eq:clt.cov}), the variance function of $\lbrace \boldsymbol{V}(t): t \geq 0 
	\rbrace$ is
		\begin{align}
			\Sigma(t) &= \text{cov}(\boldsymbol{V}(t),\boldsymbol{V}(t)) \nonumber \\
				&= \int_0^t \Phi(t,u) \cdot G(\boldsymbol{x}(u)) \cdot \Phi (t,u)^T du 
				\label{eq:clt.variance}
		\end{align}
	where $\Phi (t,s) \quad (0 \leq s \leq t)$ must satisfy equation (\ref{eq:clt.diff.phi}). 
	By differentiating equation (\ref{eq:clt.variance}) with respect to $t$,
		\begin{align}
			\frac{d \Sigma(t)}{dt} &= \Phi (t,t) \cdot G(\boldsymbol{x}(t)) \cdot \Phi (t,t)^T + 
				\frac{\partial}{\partial t} \int_0^t \Phi(t,u) \cdot G(\boldsymbol{x}(u)) \cdot \Phi 
				(t,u)^T du \nonumber \\
			&= G(\boldsymbol{x}(t)) + \partial F(\boldsymbol{x}(t)) \int_0^t \Phi (t,s) \cdot 
				G(\boldsymbol{x}(u)) \cdot \Phi (t,u)^T du + \nonumber \\
			&\qquad \int_0^t \Phi (t,u) \cdot G(\boldsymbol{x}(u)) \cdot \Phi(t,s)^T du \cdot 
				[\partial F(\boldsymbol{x}(t))]^T 
				\nonumber \\
			&= G(\boldsymbol{x}(t)) + \partial F(\boldsymbol{x}(t)) \cdot \Sigma(t) + \Sigma(t) 
				\cdot [\partial F(\boldsymbol{x}(t))]^T, \label{eq:clt.diff.variance}
		\end{align}
	where $\Sigma(0) = \boldsymbol{0}$. 

	Note that in the epidemic case, when we let $t \rightarrow \infty$, $\displaystyle \frac{d \
	Sigma(t)}{dt}$ tends to zero because if the epidemic 
	survives and takes off, it will reach the endemic stage and the infection process 
	converges to a certain distribution. 
	Otherwise, if the epidemic dies out, we will see in the later section that the stable stage is 
	reached in the disease-free stage. Both scenarios are in equilibrium level. Therefore, 
	once one of these two cases is attained, $\displaystyle \frac{d \Sigma(t)}{dt}$ definitely 
	goes to zero as $t \rightarrow \infty$.

	Now let $\hat{\Sigma}$ and $\hat{\boldsymbol{x}}$ be the variance-covariance matrix 
	and the stationary infection point when the epidemic reaches endemic stage. Hence, as 
	$t \rightarrow \infty$,
		\begin{align}
			0 = \frac{d \Sigma(t)}{dt} &= G(\hat{\boldsymbol{x}}) + \partial 
			F(\hat{\boldsymbol{x}}) \cdot \hat{\Sigma} + \hat{\Sigma} \cdot [\partial 
			F(\hat{\boldsymbol{x}})]^T \nonumber \\
			-G(\hat{\boldsymbol{x}}) &= \partial F(\hat{\boldsymbol{x}}) \cdot \hat{\Sigma} + 	
			\hat{\Sigma} \cdot [\partial F(\hat{\boldsymbol{x}})]^T. 	
			\label{eq:clt.sigma.endemic}
		\end{align}

	\begin{figure}[ht]
		\centering
		\subfloat[$(n,m, R_0) = (1000,2,5)$]{%
		\resizebox*{5cm}{!}{\includegraphics{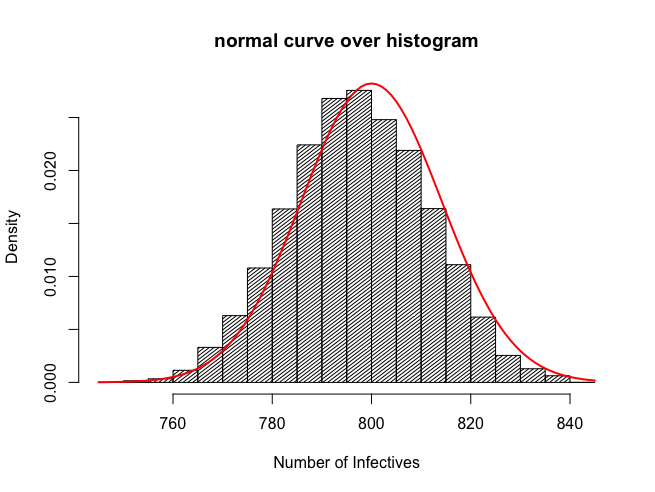}}} \hspace{5pt}
		\label{fig:sis.(n,r0)=(1e3,5)}
		\subfloat[$(n, m,R_0) = (3000,2,5)$]{%
		\resizebox*{5cm}{!}{\includegraphics{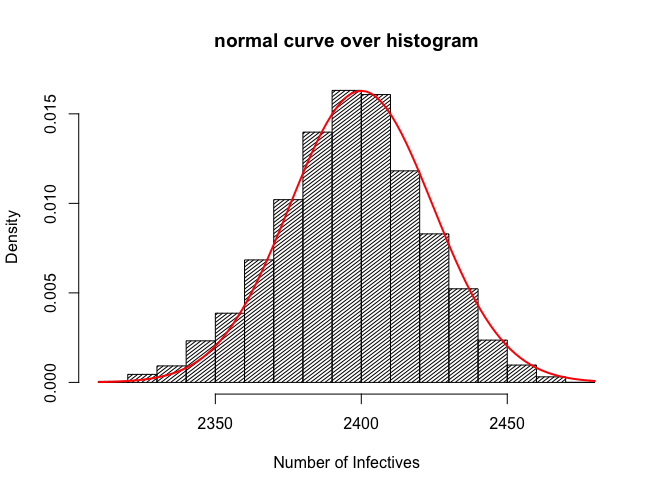}}}
		\label{fig:sis.(n,r0)=(3e3,5)} \\
		\subfloat[$(n,m, R_0) = (1000,2,8)$]{%
		\resizebox*{5cm}{!}{\includegraphics{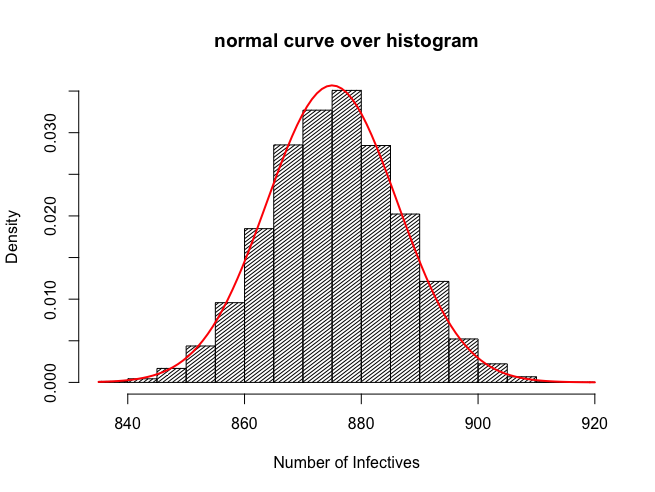}}} \hspace{5pt}
		 \label{fig:sis.(n,r0)=(1e3,8)}
		\subfloat[$(n, m,R_0) = (3000,2,8)$]{%
		\resizebox*{5cm}{!}{\includegraphics{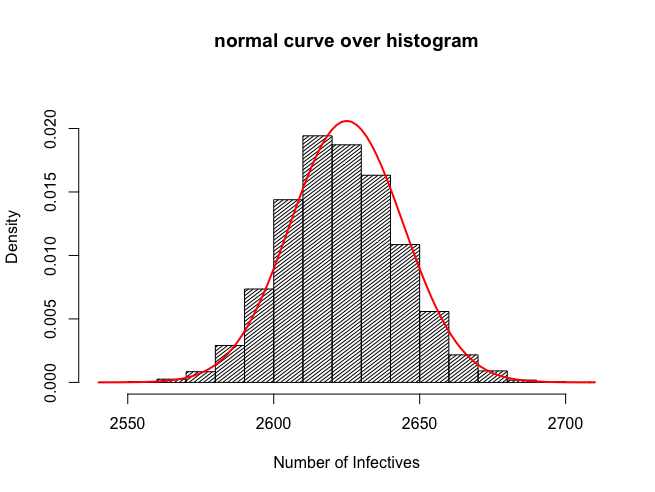}}}
		 \label{fig:sis.(n,r0)=(3e3,8)} \\
		\caption{Empirical distribution of SIS endemic level with CLT approximation}
		\label{fig:sis.empirical}
	\end{figure}

	Figure \ref{fig:sis.empirical} illustrates the distribution of $\lbrace n^{-1}Y_n(t) | 
	Y_n(t) > 0, Y_n(t) \sim \bold{\tilde{q}} \rbrace$.

\subsection{Expected Mean Time to Extinction by CLT}

	Suppose that $R_0 > 1$ and $\lbrace Y_n(t): t \geq 0 \rbrace$ is in quasi--
	equilibrium (meaning that $Y_n(t)$ enters endemic level). Then, the process will 
	definitely mimic the deterministic model and with equilibrium point  $\hat{y} = 1-
	R_0^{-1}$.

	Now, recall a drift function $F(y)$, where we define $y \in \mathbb{Z}$ in such a 
	way that $\beta_l: \mathbb{Z} \rightarrow \left[ 0,\infty \right), l \in \Delta \subseteq 
	\mathbb{Z}, \sum_l \beta_l \left( y \right) < \infty$. In SIS, $\Delta = \lbrace 
	1, -1 \rbrace$ represents the possible transitions of the process, i.e., $1$ for 
	infection to occur and $-1$ for removal to occur. Therefore, 
		\begin{align}
			F(y) &= \beta_1(y) + (-1) \cdot \beta_{-1}(y) \nonumber \\
				   &= \lambda y(t) \cdot (1-y(t)) - \gamma y(t) 
			\label{eq:sis.drift.func}
		\end{align}
	Furthermore,
		\begin{align}
			\partial F(y) = \left[ \frac{\partial F(y)}{\partial y} \right]
					&= \lambda (1-2y(t)) - \gamma \label{eq:sis.diff.drift.func}
		\end{align}
	and
		\begin{align}
			G(y) 
					&= \beta_1(y) + \beta_{-1}(y) \nonumber \\
					&= \lambda y(t) \cdot (1-y(t)) + \gamma y(t). 
					\label{eq:sis.gFunc}
		\end{align}	
	
	Thus since we assume that the process is in quasi--equilibrium, $y(t) = \hat{y}$ 
	for all $t \in \lbrace t \geq 0: Y_n(t) \mid Y_n(0) \neq 0 \sim \tilde{\boldsymbol{q}} 
	\rbrace$. Therefore,
		\begin{align}
			\partial F(\hat{y}) &=  -(\lambda - \gamma), 	
				\label{eq:sis.diff.drift.equi} \\
			G(\hat{y})  &= 2 \cdot \frac{\gamma}{\lambda} (\lambda - \gamma).
				\label{eq:sis.gFunc.equi}
		\end{align}

	Since 
	(\ref{eq:sis.diff.drift.equi}) is continuous, then  according to Theorem \ref{thm:clt}, for 
	any $t \geq 0$, 
		\begin{equation} \label{eq:sis.clt}
			\sqrt{n} \left( \frac{Y_n(t)}{n} - (1- R_0^{-1}) \right) \xrightarrow{\text{\quad D 
				\quad}} N(0, \sigma^2(t))
		\end{equation} 
	as $n \rightarrow \infty$ and $\sigma^2(t)$ must satisfy
		\begin{align}
			\frac{d \sigma^2(t)}{dt} 
			&= 2 \cdot \frac{\gamma}{\lambda} (\lambda - \gamma) - 2 \cdot (\lambda - 
				\gamma) \sigma^2(t)  \label{eq:sis.ode.sigma}
		\end{align}
	with initial value $\sigma^2(0) = 0$. Therefore, the solution of equation 
	(\ref{eq:sis.ode.sigma}), for $t \rightarrow \infty$, is
		\begin{equation}
			\sigma^2(t) = \frac{\gamma}{\lambda} \left( 1 - e^{-2(\lambda - \gamma)t} \right) 	
			\rightarrow \frac{\gamma}{\lambda} = R_0^{-1}. \label{eq:sis.sol.ode.sigma}
		\end{equation}
		
	The eq.(\ref{eq:sis.clt}) is equivalent to
		\begin{align}
			\sqrt{n} \left( \frac{Y_n(t)}{n} - (1- R_0^{-1}) \right) &\approx N(0, R_0^{-1}) 
				\nonumber \\
			\frac{Y_n(t)}{n} - (1- R_0^{-1}) &\approx N \left( 0, \frac{1}{nR_0} \right) \nonumber \\
			Y_n(t) &\approx N \left( (1- R_0^{-1}) n, \frac{n}{R_0} \right)
		\end{align}
	for large $n$ and large $t$.
	
	Suppose that $\mu_n = (1-R_0^{-1})n$ and $\sigma^2_n = \frac{\displaystyle n}
	{\displaystyle R_0}$ for simpler notations. Let $f_n$ and $F_n$ denote the pdf and 
	cdf of $N(\mu_n,\sigma_n^2)$. Then using a continuity correction, we can 
	approximate quasi-stationary distribution in (\ref{eq:sis.time.quasi}) as follows
		\begin{align}
			\tilde{q}_1 &\approx \frac{f_n(1)}{F_n \left( n + \displaystyle \frac{1}{2} \right) 
				- F_n \left( \displaystyle \frac{1}{2} \right)} \nonumber \\
			&= \displaystyle \frac{\displaystyle \frac{1}{\sqrt{2 \pi \sigma_n^2}} e^{- 
				\frac{(x-\mu_n)^2}{2 \sigma_n^2}}}{\displaystyle \int_{-\infty}^{n+\frac{1}{2}} 	
				\displaystyle \frac{1}{\sqrt{2 \pi \sigma_n^2}} e^{- \frac{(x-\mu_n)^2}{2 
				\sigma_n^2}} dx - \displaystyle \int_{-\infty}^{\frac{1}{2}} \displaystyle \frac{1}
				{\sqrt{2 \pi \sigma_n^2}} e^{- \frac{(x-\mu_n)^2}{2 \sigma_n^2}} dx} \nonumber \\
			&= \frac{ \frac{1}{\sigma_n} \phi \left( \frac{1- \mu_n}{\sigma_n} \right)}{\Phi \left( 
				\frac{n+ \frac{1}{2}-\mu_n}{\sigma_n} \right) - \Phi \left( \frac{\frac{1}{2} - \mu_n}
				{\sigma_n} \right)}
		\end{align}
	where $\phi(x)$ and $\Phi(x)$ are pdf and cdf of standardised normal distribution. 
	Therefore, the approximate mean time to extinction of the SIS epidemic model by 
	applying Theorem \ref{thm:clt} is
		\begin{align}
			E(T_Q) &= \frac{1}{\gamma \tilde{q}_1} \nonumber \\
				&\approx \left[ \frac{\displaystyle \frac{\gamma}{\sigma_n} \phi \left( \frac{ 1- 
					\mu_n}{\sigma_n} \right)}{\Phi \left( \displaystyle \frac{n+ \frac{1}{2}-\mu_n}
					{\sigma_n} \right) - \Phi \left( \displaystyle \frac{\frac{1}{2} - \mu_n}{\sigma_n} 
					\right)} \right]^{-1}
		\end{align}
	for large $n$ and large $t$.

\section{Conclusion}

	We have successfully determined the approximation of expected
	mean time to extinction of SIS model using quasi--stationary 
	distribution	. An interesting result is that the empirical distribution of 
	the SIS in endemic--stage followed normal distribution as the 
	population size went to infinity. By using the CLT of Ethier and 
	Kurtz, it turned 
	out that the mean parameters of the SIS is the 
	equilibrium points derived from solving deterministic's ODE. This result 
	was not surprising since we have showed in previous work that the 
	stochastic models moved randomly around the 
	deterministic model's equilibrium point. Furthermore, the variances of 
	the models converged to some positive constants since intuitively, as 
	the population size went sufficiently large, the stochastic paths 
	moved randomly around its endemic--equilibrium point.

\section{Appendix A} \label{app: A}

First, we need to show that for $i,j \in T$, $\boldsymbol{E}_i \boldsymbol{E}_j = 
\begin{cases}
\boldsymbol{E}_i & \text{if } i \neq j \\
\boldsymbol{0} & \text{if } i \neq j
\end{cases}. $

\begin{proof}
Note that since $\boldsymbol{u}_i$ is orthonormal, then
\begin{align}
\boldsymbol{E}_i \boldsymbol{E}_j &= \boldsymbol{W}^{-1/2} \left( \boldsymbol{u}_i  \boldsymbol{u}_i^T \right) \boldsymbol{W}^{1/2}  \boldsymbol{W}^{-1/2} \left( \boldsymbol{u}_j  \boldsymbol{u}_j^T \right) \boldsymbol{W}^{1/2} \nonumber \\
&= \left( \boldsymbol{W}^{-1/2} \boldsymbol{u}_i \right) \left(  \boldsymbol{u}_i^T \boldsymbol{u}_j \right) \left( \boldsymbol{u}_j^T \boldsymbol{W}^{1/2} \right) \nonumber \\
&= 
\begin{cases}
\boldsymbol{E}_i & \quad \text{if } i=j \\
\boldsymbol{0} & \quad \text{if } i \neq j
\end{cases} \label{eq:app.EiEj}
\end{align}
since $\boldsymbol{u}_i^T \boldsymbol{u}_j = 
\begin{cases}
1 & \quad \text{if } i=j \\
0 & \quad \text{if } i \neq j
\end{cases}. $
\end{proof}

And secondly, we need to show that $\displaystyle \sum_{i=1}^n \boldsymbol{E}_i = I_{n \times n}$.

\begin{proof}
Note that,
\begin{align}
\sum_{i=1}^n \boldsymbol{E}_i &= \left( \boldsymbol{W}^{-1/2} \boldsymbol{u}_1 \right) \left( \boldsymbol{u}_1^T \boldsymbol{W}^{1/2} \right) + \ldots + \left( \boldsymbol{W}^{-1/2} \boldsymbol{u}_n \right) \left( \boldsymbol{u}_n^T \boldsymbol{W}^{1/2} \right) \nonumber \\
&= \boldsymbol{W}^{-1/2} \left( \boldsymbol{u}_1 \boldsymbol{u}_1^T + \ldots + \boldsymbol{u}_n \boldsymbol{u}_n^T \right) \boldsymbol{W}^{1/2} \nonumber \\
&= \boldsymbol{W}^{-1/2} I_{n \times n} \boldsymbol{W}^{1/2}
= I_{n \times n}. \label{eq:app.sumE}
\end{align}
\end{proof} 

Note that $\exp{(\boldsymbol{Q}t)} = \displaystyle \sum_{k=0}^{\infty} \frac{\displaystyle \boldsymbol{Q}^k t^k}{\displaystyle k!}$. Then, for any $i \in T$,
\begin{align*}
\exp{(\lambda_i t \boldsymbol{E}_i)} &= \sum_{k=0}^{\infty} \frac{(\lambda_i t)^k \boldsymbol{E}_i^k}{k!} 
= \boldsymbol{E}_i \sum_{k=0}^{\infty} \frac{(\lambda_i t)^k}{k!} 
= \exp{(\lambda_i t)} \boldsymbol{E}_i.
\end{align*}

\section{Appendix B} \label{app: B}

Notice $\boldsymbol{E}_i = \left( \boldsymbol{W}^{-1/2} \boldsymbol{u}_i \right) \left( \boldsymbol{u}_i^T \boldsymbol{W}^{1/2} \right) $. Then,
\begin{align*}
\boldsymbol{Q}_{TT} &= \sum_{i=1}^n \lambda_i \left( \boldsymbol{W}^{-1/2} \boldsymbol{u}_i \right) \left( \boldsymbol{u}_i^T \boldsymbol{W}^{1/2} \right) \\
\left( \boldsymbol{u}_i^T \boldsymbol{W}^{1/2} \right) \boldsymbol{Q}_{TT} &= \sum_{i=1}^n \lambda_i \underbrace{ \left( \boldsymbol{u}_i^T \boldsymbol{W}^{1/2} \right) \left( \boldsymbol{W}^{-1/2} \boldsymbol{u}_i \right)}_{=1} \left( \boldsymbol{u}_i^T \boldsymbol{W}^{1/2} \right) \\
&= \sum_{i=1}^n \lambda_i \left( \boldsymbol{u}_i^T \boldsymbol{W}^{1/2} \right).
\end{align*}
Hence, $\left( \boldsymbol{u}_i^T \boldsymbol{W}^{1/2} \right)$ is the left-eigenvector of $\boldsymbol{Q}_{TT}$.

\section{References}

\end{document}